\documentclass[11pt,a4paper]{article}
\usepackage{amsmath,amscd,amssymb,amsfonts,amsthm}
\usepackage[mathscr]{eucal}
\usepackage{hyperref}

\newcommand{\semidensity}{{half-density}}
\newcommand{\semidensities}{{half-densities}}

\renewcommand{\leq}{\leqslant}

\newtheorem{thm}{Theorem}[section]

\newtheorem{lm}{Lemma}[section]
\newtheorem{prop}{Proposition}[section]
\newtheorem{cor}{Corollary}[section]

\theoremstyle{definition}
\newtheorem{de}{Definition}[section]
\newtheorem{ex}{Example}[section]

\newtheorem{rem}{Remark}[section]

\newcommand{\const}{\mathrm{const}}

\DeclareMathOperator{\Ber}{Ber} 

\DeclareMathOperator{\grad}{grad}
\renewcommand{\div}{\mathop{\mathrm{div}}}

\newcommand{\Sch}{{\mathfrak{S}}}
\newcommand{\Act}{{\mathscr{S}}}

\newcommand{\lbr}{[\![}
\newcommand{\rbr}{]\!]}
\newcommand{\bro}{{\boldsymbol{\rho}}}
\newcommand{\brro}{{{\boldsymbol{\rho}}{\mathrm{'}}}}

\DeclareMathOperator{\ad}{ad}
\newcommand{\lie}[1]{{\cal L}_{{#1}}}
\newcommand{\liex}[1]{{\lie{{#1}}}}
\newcommand{\DD}[1]{\Delta_{#1}}
\DeclareMathOperator{\Dr}{\DD{{\bro}}}
\DeclareMathOperator{\DDr}{\DD{{\bro}}^2}
\newcommand{\Drr}{\DD{\brro}\,}
\newcommand{\DDrr}{\DD{\brro}^2\,}

\DeclareMathOperator{\divr}{{\mathrm{div}_{\bro}}}

\newcommand{\der}[2]{{\frac{\partial {#1}}{\partial {#2}}}}

\newcommand{\dder}[3]{{\frac{\partial^2 {#1}}{\partial {#2}\partial {#3}}}}

\newcommand{\fun}{C^{\infty}}

\def\e{\varepsilon}
\def\s{\sigma}
\def\f{{\varphi}}

\def\D{\Delta}

\renewcommand{\r}{{\rho}}

\def\d{\delta}
\def\t{\theta}
\def\l{\lambda}

\newcommand{\ft}{{\tilde f}}
\newcommand{\gt}{{\tilde g}}

\newcommand{\at}{{\tilde a}}
\newcommand{\bt}{{\tilde b}}

\newcommand{\Xt}{{\tilde X}}

\newcommand{\Pt}{{\tilde P}}

\newcommand{\bs}{{\boldsymbol{s}}}
\newcommand{\bu}{{\boldsymbol{u}}}

\numberwithin{equation}{section}

\title{On odd Laplace operators
}

\newcommand{\myaddress}{{\vspace{-0.4cm}\footnotesize
$^{\scriptstyle{1}}$ Department of Mathematics, University of
Manchester Institute of Science and Technology (UMIST), PO Box 88,
Manchester M60 1QD, England\\ \smallskip} {\footnotesize
$^{\scriptstyle{2}}$ G.~S.~Sahakian~Department~of~Theoretical
~Physics, Yerevan State University, \\1 A. Manoukian Street,
375049 Yerevan, Armenia\\} {\footnotesize {\tt
theodore.voronov@umist.ac.uk, khudian@umist.ac.uk} }}

\author{Hovhannes M. Khudaverdian{\small $^{\scriptstyle{1,2}}$},
Theodore  Voronov{\small $^{\scriptstyle{\,1}}$}}
\date{\myaddress}

\begin{document}
\maketitle \vspace{-1.1cm}
\begin{abstract}\footnotesize
We consider odd Laplace operators acting on densities of various
weights on an odd Poisson (= Schouten) manifold $M$. We prove that
the case of densities of weight $1/2$ (\semidensities) is
distinguished by the existence of a unique odd Laplace operator
depending only on a point of an ``orbit space'' of volume forms.
This includes earlier results for the odd symplectic case, where
there is a canonical odd Laplacian on \semidensities. The space of
volume forms on $M$ is partitioned into orbits by the action of a
natural groupoid whose arrows  correspond to the solutions of the 
quantum Batalin--Vilkovisky equations.  We compare this situation
with that of Riemannian and even Poisson manifolds. In particular,
we show that the square of an odd Laplace operator is a Poisson
vector field defining an analog of Weinstein's ``modular class''.
\end{abstract}

\vspace{-0.9cm} \tableofcontents
\section*{Introduction}
In this paper we study odd Laplace operators acting on densities
of various weights on odd Poisson manifolds. This happens to be a
very rich geometrical topic, surprisingly linking even and odd
Poisson geometry with Riemannian geometry. The main result is the
construction of odd Laplacians acting on \semidensities { } (=
densities of weight $1/2$) in the context of Poisson geometry and
a detailed study of their properties.  

Earlier in~\cite{hov:max, hov:semi} it was shown that a canonical
odd Laplacian on \semidensities { } exists on odd symplectic
manifolds. In the current paper we consider Laplacians on
\semidensities { }on arbitrary odd Poisson manifolds. This change
of viewpoint has led us to a completely new picture.

Namely, rather than study a single Laplace operator, we consider
all operators depending on arbitrary volume elements on the
supermanifold. Laplace operators acting on \semidensities{}
are characterized  by a commutator identity of the form $[\D,f]=
\liex{f}$, where $\liex{f}$ is the Lie derivative along the
Hamiltonian vector field with  Hamiltonian $f$. Starting from this
relation, we arrive at a natural groupoid whose orbits parametrize
the set of odd Laplace operators: \semidensities { }are distinguished
because for them the Laplacian depends not on a volume element but only
on its orbit. This groupoid, which we call the ``master groupoid'',
appears also in other instances. It seems that it plays a very
important role in odd Poisson geometry. Another interesting new
object is a remarkable analog of Weinstein's ``modular class'',
which comes from the square of odd Laplace operators.

Besides naturally expected links with even Poisson geometry, we
have found far-reaching surprising links with Riemannian geometry. We
found it worthy to point out at a simple analogy between the
Batalin--Vilkovisky formalism in quantum field theory and the
usual quantum mechanics, in particular, between the exponential of
the ``quantum master action'' and  the \semidensity { }wave
function.

\textit{Remark.} The Batalin--Vilkovisky formalism appeared around
1981 as a very general method of quantization of systems with
gauge freedom~\cite{bv:perv, bv:vtor, bv:closure}. The central
role in this method is played by a second order odd differential
operator somewhat similar to the divergence of multivector fields.
From the viewpoint  of supermanifolds  it is better to interpret
it as an analog of a Laplacian. Geometry of the
Batalin--Vilkovisky (BV) formalism has been investigated in many
works, in particular~\cite{ass:bv} and \cite{hov:deltabest,
hov:khn1}. Algebraic structures related with the BV formalism
gradually became very fashionable. Invariant geometric
constructions for the BV operator have been studied
in~\cite{hov:deltabest}, \cite{ass:bv} in the symplectic case and
recently in \cite{yvette:divergence} in the Poisson case. In these
works the operators considered act on functions. A new approach 
based on \semidensities { }was suggested in~\cite{hov:max,
hov:semi} and is developed here.

\section{Some facts from odd Poisson geometry}
\subsection{Poisson brackets and Hamiltonian vector fields}

Let us briefly recall some well known formulae, for  reference
purposes\footnote{A  collection of useful formulae can be found in
the introductory section of~\cite{tv:graded}}. Given a
supermanifold $M$, an arbitrary \textit{odd Poisson} (or
\textit{Schouten}) \textit{structure} on $M$ is specified by an
odd function on $T^*M$ quadratic in momenta:
\begin{equation}\label{oddbracket}
    \Sch=\frac{1}{2}\,\Sch^{ab}(x)p_bp_a
\end{equation}
(``odd quadratic Hamiltonian''). Explicitly:
\begin{equation}\label{eqfgS}
    \{f,g\}_{\Sch}:=(f,(\Sch,g))=((f,\Sch),g)
    =-(-1)^{\ft(\at+1)}\Sch^{ab}\der{f}{x^b}\der{g}{x^a},
\end{equation}
where we denote by $(\ ,\ )$ the canonical even Poisson bracket on
$T^*M$. Here $\{f,g\}_{\Sch}$ stands for the odd bracket on $M$
specified by $\Sch$. The Jacobi identity for $\{f,g\}_{\Sch}$ is
equivalent to the vanishing of the canonical Poisson bracket
$(\Sch,\Sch)=0$   on $T^*M$.

In the sequel we leave $\Sch$ fixed and drop the reference to
$\Sch$ from the notation for the bracket. The Hamiltonian vector
field on $M$ corresponding to a function $f$ is defined as
\begin{equation}\label{eqxf}
    X_f:=(-1)^{\ft+1}\{f,\ \}=(-1)^{\at\ft}\Sch^{ab}\der{f}{x^b}\der{}{x^a}.
\end{equation}
It has parity opposite to that of $f$. Notice that
\begin{equation}\label{eqxskobka}
    X_{\{f,g\}}=[X_f,X_g]
\end{equation}
and
\begin{equation}\label{xfg}
    X_{fg}=(-1)^{\ft}fX_g+(-1)^{\gt+\ft\gt}gX_f
\end{equation}
(or    $X_{fg}=X_f\circ g+(-1)^{\ft}fX_g+(-1)^{\ft}\{f,g\}$, where
in the r.h.s. there is the multiplication by $g$ followed by the
action of $X_f$).

We denote by
\begin{equation}
    \liex{f}:=\lie{X_f}
\end{equation}
the Lie derivative of geometric objects on $M$ (e.g., tensor
fields of a given type) w.r.t. the vector field  $X_f$. The Lie
derivative along $X_f$ makes the space of particular geometric
objects on $M$  into a ``Poisson module'' over the odd Poisson
algebra $\fun(M)$, meaning that
\begin{align}
    [\liex{f},g]&=(-1)^{\ft+1}\{f,g\}=X_fg\\
    [\liex{f},\liex{g}]&=\liex{\{f,g\}}.
\end{align}

\subsection{Odd Laplacian on functions}

Select a volume form $\bro=\r\, Dx$ on an odd Poisson manifold
$M$. We assume that $\bro$ is non-degenerate, i.e., that it
constitutes a basis element of the space of volume forms. (Here    
and in the sequel we skip the question of global existence of 
such a form  
and the questions of orientation. It is not difficult to supplement the corresponding details.) 

\begin{de}
The \textit{Laplace operator} acting on functions  is defined by
the formula
\begin{equation}\label{eqlaplacefunc}
    \Dr f:=\text{``$\div\grad f$''}=\divr X_f=
    \frac{\liex{f}\bro}{\bro}.
\end{equation}
It depends on the
choice of $\bro$. The operator $\Dr$ is odd.
\end{de}

Formula~\eqref{eqlaplacefunc} was introduced for the first time in
1989 in~\cite{hov:deltabest} (formally, only for the symplectic 
case) as an invariant construction for the ``$\D$-operator'' of
Batalin and Vilkovisky~\cite{bv:perv,bv:closure}\footnote{As we
see it now, the physical meaning requires rather an operator
acting on \semidensities, i.e., densities of weight $1/2$. Such
operators are the main topic of this paper --- though we start off
from operators on functions.}. See 
also~\cite{hov:khn1},\cite{ass:bv}. A detailed analysis of this  
construction for the Poisson case, in a very 
general algebraic setup,  is given in~\cite{yvette:divergence}. 

From formula~\eqref{eqxf} for Hamiltonian vector fields and the
formula for the divergence
\begin{equation*}
    \divr X=\frac{1}{\r}\,(-1)^{\at(\Xt+1)}\der{(\r X^a)}{x^a}
\end{equation*}
one immediately gets a simple ``Laplace--Beltrami type''
expression for $\Dr$:
\begin{equation}\label{eqdeltacoor}
    \Dr f=\frac{1}{\r}\,\der{}{x^a}\left(\r\,
    \Sch^{ab}\der{f}{x^b}\right)=\Sch^{ab}\dder{f}{x^b}{x^a}+
    \text{(lower order terms)}.
\end{equation}
The operator $-\frac{\hbar^2}{2}\Dr$ is a ``quantization'' of the
function $\Sch$ on $T^*M$. 
Notice that there is a coordinate-free expression similar
to~\eqref{eqfgS}
\begin{equation}
    \Dr f\cdot \f_{\bro}=((\Sch,f),\f_{\bro})
\end{equation}
via the canonical \textit{even} Poisson brackets on $T^*M$, where
$\f_{\bro}(x,p)=\r(x)\d(p)$ is the generalized function on $T^*M$
corresponding to the volume form $\bro$ on $M$. Related to it is
another useful coordinate-free representation of $\Dr$  given by
an integral identity
\begin{equation}\label{eqgreen}
    \int_M f(\Dr g)\,\bro=\int_M \{f,g\}\,\bro,
\end{equation}
analogous to the familiar ``Green's first integral formula'' in
Riemannian geometry (without boundary terms). It is valid if $f$
or $g$ is compactly supported inside $M$. It follows that $\Dr$ is
formally self-adjoint with respect to $\bro$.

{}From the definition of $\Dr$, together with the
formula~\eqref{eqxskobka}, one easily obtains the derivation
property with respect to the bracket:
\begin{equation}\label{eqlapskobka}
    \Dr \{f,g\}= \{\Dr f,g\}+(-1)^{\ft+1}\{f,\Dr g\}.
\end{equation}
Using the definition of $X_f$, this can be rewritten in the
commutator form
\begin{equation}\label{eqkom}
    [\Dr,X_f]=-X_{\Dr f}.
\end{equation}
As for the product of functions, the odd Laplace operator
$\D_{\bro}$ satisfies 
\begin{equation}\label{eqkomdliafun1}
    \Dr (fg)= (\Dr f)g+(-1)^{\ft}f(\Dr g)+(-1)^{\ft+1}2\{f,g\}
\end{equation}
(directly analogous to  another ``Green's identity''). It follows
from formula~\eqref{xfg} and the equality   $\divr(fX)=f\divr
X+(-1)^{\Xt\ft}Xf$. Recalling the definition of $X_f$, we
rewrite~\eqref{eqkomdliafun1} in the commutator form
\begin{equation}\label{eqkomdliafun2}
    [\D_{\bro},f]=2X_{f}+\D_{\bro}f.
\end{equation}
From here  purely algebraically follow  the identities
\begin{gather}
    \Dr (f^n)=nf^{n-1}\Dr\!f-n(n-1)f^{n-2}\{f,f\}, \\
    \Dr e^{kf}=k\Bigl(\Dr f-k\,\{f,f\}\Bigr)\,e^{kf}
    \label{eqlapexp}.
\end{gather}

The definition of the Laplace operator $\Dr$ depends on the choice
of a basis volume form $\bro$. If we change $\bro$, the odd
Laplacian on functions transforms as
\begin{equation}\label{eqzamenafun}
    \Drr = \Dr +X_{\s}=\Dr-\{\s,\ \},
\end{equation}
where $\brro=e^{\s}\bro$.

\subsection{Symplectic case: a canonical  odd Laplacian on \semidensities}

The case of a nondegenerate bracket, i.e., of an \textit{odd
symplectic manifold} $M$, is distinguished by the existence of
local Darboux coordinates, i.e., local coordinates $x^i,\t_i$
($x^i$ even, $\t_i$ odd) such that the odd bracket has the form
$\{\t_i,x^j\}=-\{x^j,\t_i\}=\d_i^j$ and
$\{\t_i,\t_j\}=\{x^i,x^j\}=0$.
Suppose 
$\bro=\r(x,\t)\,{D(x,\t)}$ in such coordinates. Then, clearly,
\begin{gather}
    \Dr=\Delta_0+X_{\ln\r}=\Delta_0-\{\ln\r,\ \}, \label{eqcoordeltarho}\\
    \intertext{where}
    \Delta_0=2\dder{}{x^i}{\t_i} \label{eqcoordeltarho2}
\end{gather}
is a ``coordinate $\D$-operator''.  Each term
in~\eqref{eqcoordeltarho} is not invariant under changes of
coordinates, only the sum being invariant. It is tempting,
however, to consider $\D_0$ independently. It turns out that
though it  does not make sense as an operator acting on functions,
an analog of $\D_0$ is well-defined as a canonical operator on
densities of weight $1/2$ (\semidensities).

\begin{de}[\textnormal{see~\cite{hov:max,hov:semi}}]
Let $\bs$ be a \semidensity.  The \textit{canonical odd Laplace
operator} $\D$ acts on $\bs$ by the formula
\begin{equation}\label{eqcanlap}
    \D \bs:=\left(2\dder{s}{x^i}{\t_i}\right){D(x,\t)}^{1/2},
\end{equation}
where $\bs=s(x,\t)\,{D(x,\t)}^{1/2}$ in a local Darboux chart.
\end{de}

The operator~\eqref{eqcanlap} was introduced in~\cite{hov:max}. It
was proved in~\cite{hov:max,hov:semi} that the definition of $\D$
on \semidensities { }does not depend on the choice of a Darboux
chart, thus yielding a well-defined operator, canonical in the
sense that it depends only on the odd symplectic structure and
does not require any extra data like a volume form (in contrast to
the operator $\Dr$ on functions). The role of \semidensities{ } is
crucial. One can show that an operator defined in Darboux
coordinates by a formula similar to~\eqref{eqcanlap} on arbitrary
densities of weight $w$ will be invariant \textit{only} for
$w=1/2$.

The existence of $\D$ on half-densities is
essentially equivalent to the following statement, which can be
traced  back to Batalin and Vilkovisky~\cite{bv:closure}:

\begin{lm}[``Batalin--Vilkovisky Lemma''] \label{thmbvlemma}
\begin{equation}\label{eqbvlemma}
    \D_{0}\left(\Ber{\der{x'}{x}}\right)^{1/2}=0
\end{equation}
for the change of coordinates between two Darboux charts. Here
$x=(x^i,\t_i)$, $x'=(x^{i'},\t_{i'})$, and  $\D_0$ corresponds to
the ``old'' coordinate system.
\end{lm}

We do not give a proof here. Proofs and a detailed analysis of the
properties of the operator $\D$ on \semidensities { }can be found
in~\cite{hov:semi}. Further analysis will be given in
subsection~\ref{subsecanal}.

The exceptional role of the exponent $1/2$ in
equation~\eqref{eqbvlemma}  cannot be detected on the
infinitesimal level. It is related with the possibility to
``integrate'' infinitesimal canonical transformations to finite
ones, due to a deep \textit{groupoid property} of the
Batalin--Vilkovisky equations we discuss in Section 2.

The following important formulae were obtained in the calculus of
\semidensities { }on an odd symplectic manifold~\cite{hov:semi}:

\begin{equation}\label{eqlapnafunplot}
    \D(f\bs)=(\D_{\bs^2}f)\bs + (-1)^{\ft}f(\D\bs)
\end{equation}
for an arbitrary function $f$ and a non-degenerate \semidensity {
}$\bs$ (at the r.h.s. stands the Laplacian on functions with respect to
the volume form $\bs^2$), and
\begin{equation}\label{eqdkvadratnafunk}
    \D^2_{\bro}f=\left\{{\bro}^{-1/2}\,\D({\bro}^{1/2}),f\right\},
\end{equation}
where on the l.h.s. stands the Laplace operator on functions
with respect to the volume form $\bro$ and at the r.h.s. stands
the canonical operator on \semidensities.
Equation~\eqref{eqdkvadratnafunk} in particular implies that the
square of the odd Laplacian $\D_{\bro}$ is a Hamiltonian vector
field. Equation~\eqref{eqlapnafunplot} can be restated as follows:
\begin{equation}\label{eqlapfunplot2}
    \D(f\bs)=2\liex{f}\,\bs + (-1)^{\ft}f\,\D\bs,
\end{equation}
now valid for arbitrary \semidensities { }(indeed, for an even
half-density $\bs$,
$(\D_{\bs^2}f)\bs=(\liex{f}\bs^2)\bs^{-1}=2\liex{f}\,\bs$). Notice
that equation~\eqref{eqlapfunplot2} means that
\begin{equation}\label{eqlapfunplot3}
    [\D,f]=2\liex{f},
\end{equation}
(compare with equation~\eqref{eqkomdliafun2} for functions).

\medskip
\textit{Physical background.}  By the geometrical meaning of the
Batalin--Vilko\-visky quantization procedure, the exponential of
the ``quantum master action'' $e^{i\Act/\hbar}$ appearing in the
quantum master equation $\D\,e^{i\Act/\hbar}=0$ is not a scalar,
but the coefficient of a density of weight $1/2$ on the extended
phase space of fields, ghosts, antifields, antighosts. The
$\D$-operator in this master equation exactly corresponds to the
canonical Laplace operator on \semidensities.  In the BV-method
$e^{i\Act/\hbar}$ is integrated over a Lagrangian submanifold;
that is due to the fact that  \semidensities { }on the phase space
correspond to forms on Lagrangian submanifolds~\cite{hov:semi}. The
difference of the BV quantum action (logarithm of a \semidensity)
from a scalar appears in quantum corrections to the scalar
classical action. In the usual Feynman integral (without gauge
freedom) the exponential of the ``quantum action'' $e^{iS/\hbar}$
is the coefficient of a volume form, i.e., of a density of weight $1$.

\section{Odd Laplace operators: main results}\label{secmain}

\subsection{Square of $\Dr$,  master groupoid and
 modular class} \label{subsecmaster}

In the \textit{even} Poisson case, the analog of
formula~\eqref{eqlaplacefunc} gives an operator of the first
order, a Poisson vector field (see~\cite{weinstein:modular} and
references therein).  The reason why there are no terms of the
second order is the skew-symmetry of the Poisson tensor, in
contrast with the symmetry of $\Sch^{ab}$. In fact, it is nothing
but the divergence of the Poisson bivector with respect to a chosen
volume form.  Under a change of the volume form, this Poisson vector
field changes by a Hamiltonian vector field. Hence, its class in the
Poisson--Lichnerowicz cohomology does not depend on the volume form
and is an invariant of even Poisson manifolds (Weinstein's
``modular class'').

Now we shall see how a similar class arises for odd Poisson
(Schouten) manifolds.

Consider the square of the Laplace operator $\D_{\bro}$. Applying
$\Dr$ to both sides of~\eqref{eqkomdliafun1} and using the
derivation property with respect to  the bracket~\eqref{eqlapskobka}, after
cancellations we arrive at the simple formula
\begin{equation}\label{eqkvadrat}
    \DDr(fg)=(\DDr f)g+(-1)^{\ft} f(\DDr g),
\end{equation}
which shows that $\DDr$ is a vector field. Notice that $\Dr$ is 
of  order $\leq 2$, so $\DDr$ might, in principle,
contain terms of order $\leq 4$. However, since
$\DDr=(1/2)[\Dr,\Dr]$, the terms of order $4$ cancel
automatically. The cancellation of the terms of order $3$ follows from
the vanishing of the canonical bracket $(\Sch,\Sch)$, which is the
classical limit of $[\Dr,\Dr]$. These simple considerations
\textit{a priori} allow to reduce the order of $\DDr$ to $2$.
Remarkably, the actual order is $1$.

(In an algebraic setup, the equivalence of $\D$ being a
derivation of the odd bracket and $\D^2$ being a derivation of
the associative product  was explicitly noticed
in~\cite{yvette:divergence}.)

Since $\DDr=(1/2)[\Dr,\Dr]$ and  $\Dr$ is a derivation of the odd
bracket, we arrive at
\begin{prop}
The vector field $\DDr$ is  a derivation of the odd bracket, i.e.,
$\DDr$ is a Poisson vector field.
\end{prop}

Recall that in the symplectic case it is always Hamiltonian
(equation~\eqref{eqdkvadratnafunk}).

\begin{de}
The Poisson vector field $\DDr$ will be called  \textit{the
modular field} of the Schouten manifold $M$ with respect to the
volume form $\bro$.
\end{de}

\begin{prop}
The modular field $\DDr$ preserves the volume form $\bro$.
\end{prop}
\begin{proof}
Denote $Y:=\DDr$. To calculate the Lie derivative $\lie{Y}\bro$,
we can apply ``Green's formula''~\eqref{eqgreen} (or the
self-adjointness of $\Dr$):
\begin{equation*}
    \int_M g(\lie{Y}\bro)=-\int_M (Yg)\bro=-\int_M (\DDr
    g)\bro=-\int_M (\Dr(1)\Dr g)\bro=0,
\end{equation*}
for an arbitrary test function $g$. Hence $\lie{Y}\bro=0$.
\end{proof}

Consider now the transformation of $\DDr$ under a change of $\bro$.

\begin{thm}\label{thmzamenakvadrat}
If $\bro'=e^{\s}\bro$, then
\begin{equation}\label{eqzamenakvadrat}
    \DDrr =\DDr- X_{H(\bro',\bro)},
\end{equation}
where
\begin{equation}\label{eqhproro}
    H(\bro',\bro):=\Dr {\s}-\frac{1}{2}\,\{{\s},{\s}\}=2e^{-{{\s}}/{2}}\Dr
    e^{{{\s}}/{2}}.
\end{equation}
The odd functions $H(\bro',\bro)$ satisfy the following ``cocycle
conditions'':
\begin{gather}
    H(\bro,\bro)=0, \label{eqhroro}\\
    H(\bro',\bro)+H(\bro,\bro')=0, \label{eqhkos}\\
    H(\bro'',\bro)=H(\bro'',\bro')+H(\bro',\bro). \label{eqhkoc}
\end{gather}
\end{thm}
\begin{proof}
By equation~\eqref{eqzamenafun}, we have $\Drr=\Dr+X_{\s}$. Hence
\begin{multline*}
    \DDrr=(\Dr+X_{\s})^2=\DDr+[\Dr,X_{\s}]+X_{\s}^2=
    \DDr-X_{\Dr {\s}}+\frac{1}{2}\,[X_{\s},X_{\s}]=\\
    \DDr-X_{\Dr {\s}-\frac{1}{2}\,\{{\s},{\s}\}},
\end{multline*}
which proves~(\ref{eqzamenakvadrat}-\ref{eqhproro}). The cocycle
conditions are checked directly:
\begin{multline*}
    H(\bro'',\bro')+H(\bro',\bro)=
    \Drr \s'-\frac{1}{2}\,\{\s',\s'\}+\Dr
    {\s}-\frac{1}{2}\,\{{\s},{\s}\}=\\
    \Dr \s'-\{\s,\s'\}-\frac{1}{2}\,\{\s',\s'\}+\Dr {\s}-
    \frac{1}{2}\,\{{\s},{\s}\}=\\
    \Dr (\s+\s')-\frac{1}{2}\,\{\s+\s',\s+\s'\}=H(\bro'',\bro).
\end{multline*}
Similarly  $H(\bro,\bro')=-H(\bro',\bro)$, and $H(\bro,\bro)=0$ is
obvious.
\end{proof}

\begin{rem}
The vector fields $X_{H(\bro',\bro)}=\DDr-\DDrr$ being
``coboundaries''  trivially satisfy the ``cocycle condition''.
However, one might expect that the corresponding Hamiltonians
$H(\bro',\bro)$ define a cocycle  only up to Casimirs; the key
statement is that the cocycle conditions are satisfied exactly.
\end{rem}

Consider the factors $e^{\s}$ as ``arrows'' between volume forms
on $M$. This is an action of a group on the set of volume forms.
Consider now only the arrows $\bro\to\bro'=e^{\s}\bro$ satisfying
the equation $\Dr e^{\s/2}=0$ (the equation on an arrow depends on
its ``source''). From formula~\eqref{eqhproro} and the cocycle
conditions~(\ref{eqhroro}), (\ref{eqhkos}), (\ref{eqhkoc}) follows
an important statement.

\begin{thm}\label{thmgrupoid}
The  solutions of the equations $\Dr e^{\s/2}=0$ form a groupoid.
That is: for three volume forms $\bro,\bro',\bro''$, if $\Dr
e^{\s/2}=0$ and $\Drr e^{\tau/2}=0$, then $\Dr e^{(\s+\tau)/2}=0$
(``composition''). Also, if $\Dr e^{\s/2}=0$, then $\Drr
e^{-\s/2}=0$ (``inverses''). Here $\bro'=e^{\s}\bro$,
$\bro''=e^{\tau}\bro'$.
\end{thm}

Since the arrows in it satisfy equations that have the form of the
quantum ``master equation'' of the Batalin--Vilkovisky formalism,
we shall call the groupoid defined in Theorem~\ref{thmgrupoid} the
\textit{master groupoid} of an odd Poisson manifold $M$.  The
space of all non-degenerate volume forms on $M$ is partitioned
into orbits of the master groupoid.

The factor $1/2$ in the exponent is exceptional; no such property
holds for any $e^{\l\s}$ other than $\l=1/2$.

\begin{ex}
Let $M$ be an odd symplectic manifold. Consider as points of a
groupoid all Darboux coordinate systems and as arrows the
canonical transformations between them. A homomorphic image of it
is the groupoid whose points are coordinate volume forms (in
Darboux coordinates) and whose arrows are the respective Jacobians
$J$. Consider a new groupoid with arrows $J^{\l}$.
Infinitesimally, $\D_0 J^{\l}=0$ for any $\l$, where
$J=1+\e\div_{0}X_{F}=1+\e\D_0F$,  where $F$ is a Hamiltonian  
generating the canonical transformation, because $\D_0^2=0$. 
However, we can glue together the conditions for infinitesimal
transformations $\D_0J^{\l}=0$ only when $\l=1/2$. Thus we arrive
at the identity $\D_0 J^{1/2}=0$ for a finite transformation,
i.e., to the Batalin--Vilkovisky Lemma. That means that all
``Darboux coordinate volume forms'' belong to the same orbit of
the master groupoid, i.e., in the symplectic case the orbit space
has a natural base point.
\end{ex}

Notice now that by Theorem~\ref{thmzamenakvadrat} the  modular
vector field $\D_{\bro}^2$ depends on a volume form up to a
Hamiltonian vector field. Hence it defines a cohomology class
$[\D_{\bro}^2]$ depending only on the odd bracket structure. We
call $[\D_{\bro}^2]$ the \textit{modular class} of the Schouten
manifold $M$. More precisely, consider $T^*M$ with the canonical
even bracket. The Schouten tensor~\eqref{oddbracket} satisfies
$(\Sch,\Sch)=0$. Hence the operator $D:=\ad \Sch=(\Sch,\ )$ on the
space $\fun(T^*M)$  is an odd differential. We  call the complex
$(\fun(T^*M),D)$ or its subcomplex consisting of fiberwise
polynomial functions, the \textit{Schouten--Lichnerowicz complex}
of $M$, and its cohomology, the \textit{Schouten--Lichnerowicz
cohomology} of $M$. The modular class belongs to the first
cohomology group.  For odd symplectic manifolds (constant Schouten
structure) this class vanishes by
equation~\eqref{eqdkvadratnafunk}. One can show that it also
vanishes for all linear Schouten structures. (Compare with
Weinstein's class for the Berezin bracket on $\mathfrak{g}^*$ if
$\mathfrak{g}$ is non-unimodular.) Whether there are Schouten
structures with a nontrivial modular class, is an open question.

\subsection{Laplacians on \semidensities{} (Poisson case)}

In the symplectic case the canonical odd Laplacian acting on
\semidensities { }satisfies equation~\eqref{eqlapfunplot3}. We
shall use this identity  to characterize Laplacians acting on
\semidensities { }on an arbitrary odd Poisson manifold.

\begin{thm}\label{thmdonsemi} 
Given an odd Poisson manifold, a linear differential operator
$\D$ acting on densities of weight $w$ and satisfying  the
equation
\begin{equation}\label{eqkomdliapolu}
    [\D,f]=2\liex{f}
\end{equation}
for an arbitrary function $f$, exists if and only if $w=1/2$. In
this case $\D$ is defined uniquely up to a zeroth-order term. If
$\bs=s\bro^{1/2}$, where $\bro$ is some basis volume form, then
\begin{equation}\label{eqlapnapolu1}
    \D(\bs)=\D(s\bro^{1/2})=(\D_{\bro}s)\bro^{1/2}+C\bs
\end{equation}
(where $C$ is a function). Here $\D_{\bro}$ is the Laplace
operator  on functions.
\end{thm}
\begin{proof}
Consider densities of weight $w$ and suppose an operator with the
property~\eqref{eqkomdliapolu} exists. For any two operators
satisfying~\eqref{eqkomdliapolu} their difference commutes with
the multiplication by functions, hence is a scalar. To fix this
scalar, we can  set $\D(\bro^w)=0$ for some chosen volume form
$\bro$. Then from~\eqref{eqkomdliapolu} immediately follows that
$\D(\bs)=\D(s\bro^{w})=2\liex{s}\bro^w=
2w\liex{s}\bro\cdot\bro^{w-1}=2w(\D_{\bro}s)\bro\bro^{w-1}=
2w(\D_{\bro}s)\bro^{w}$. (In particular, for $w=1/2$ we get
exactly $\D(\bs)=(\D_{\bro}s)\bro^{w}$.) However,  a direct check
shows that the operator defined by this formula satisfies
condition~\eqref{eqkomdliapolu} only for $w=1/2$.
For all other weights the actual commutator contains an extra term.
We omit here these calculations.
\end{proof}

Choose a volume form $\bro$ and  fix normalization by requiring
$\D(\bro^{1/2})=0$.  We arrive at the following definition.
\begin{de}\label{deflapnapolu}
The \textit{odd Laplace operator on \semidensities { }on an odd
Poisson manifold} is
\begin{equation}\label{eqlapnapolu2}
    \D(\bs):=\bro^{1/2}\D_{\bro}(\bs\bro^{-1/2})
\end{equation}
where $\bro$ is some basis volume form.
\end{de}

Thus defined, the Laplace operator $\D$ on \semidensities {
}satisfies the commutator
condition~\eqref{eqkomdliapolu}\footnote{By multiplying $\D$ by a
constant the factor ``$2$'' in front of the Lie derivative
in~\eqref{eqkomdliapolu} can be replaced by  any number. It is
convenient, though, to keep coherent normalization with the
Laplacians on functions.}.  The uniqueness implies that varying of
$\bro$ in~\eqref{eqlapnapolu2} changes $\D$ by a scalar term, in
contrast with \eqref{eqzamenafun} for the Laplacian on functions.

\begin{thm}\label{thmzamenanapolu}
Under the change of volume form $\bro \to \bro'=e^{\s} \bro$, the
Laplace operator on \semidensities { }transforms as follows:
\begin{equation}\label{eqzamenarodliapolu1}
    \D' =\D  -e^{-{{\s}}/{2}}\Dr e^{{{\s}}/{2}} =\D
    -\frac{1}{2}\,H(\bro',\bro),
\end{equation}
where we denote by $\D'$ the operator corresponding to the volume
form $\bro'$. Here $H(\bro',\bro)$  is the Hamiltonian \ref{eqhproro}.
\textnormal{(Notice that the formula
$\D'-\D=-\frac{1}{2}\,H(\bro',\bro)$ immediately proves that
$H(\bro',\bro)$ is a cocycle.)}
\end{thm}
\begin{proof}
From Theorem~\ref{thmdonsemi} follows that the operator $\D'-\D$
commutes  with all functions. Thus $\D'-\D=V(\bro',\bro)$ must be
a function depending on $\bro and \bro'$. To find it, apply
$\D'-\D$ to ${\bro'}^{1/2}$ and use the normalization condition
$\D'({\bro'}^{1/2})=0$. We obtain that
$-\D({\bro'}^{1/2})=V(\bro',\bro)\,{\bro'}^{1/2}$, i.e., by the
definition of the Laplacian on \semidensities,
$-(\D_{\bro}e^{\s/2})\,\bro^{1/2}=V(\bro',\bro)\,e^{\s/2}\bro^{1/2}$,
hence    $V(\bro',\bro)=-e^{-\s/2}\D_{\bro}e^{\s/2}$.
\end{proof}

Recall that the condition $\Dr e^{{{\s}}/{2}}=0$
specifies the arrows of the master groupoid.

\begin{cor}
The odd Laplace operator on \semidensities { }is constant on the
orbits of  the master groupoid.
\end{cor}

The operator $\D$ actually depends not on a volume form $\bro$,
but only on its orbit under the action of the master groupoid. The
situation is drastically different from that for Laplacians on
functions and for densities of weight $\neq 1/2$.

The master groupoid, which appeared above in relation with the
transformation law of $\D_{\bro}^2$ on functions, now directly
arises from the transformation law of $\D$ on \semidensities. In
terms of operators on \semidensities, the defining equation
$\Dr e^{\s/2}=0$ takes the transparent form $\D({\bro'}^{1/2})=0$.

\subsection{Analysis of the symplectic case}\label{subsecanal}

Now we can briefly review the symplectic case. It is distinguished
by the existence of Darboux charts. With every such chart one can
associate a coordinate volume form $\bro=D(x,\t)$. Moreover, one 
can construct a global volume form such that for some Darboux  
atlas 
this form coincides with the 
coordinate volume form in every chart of this atlas. The  
proof uses some facts about the topology of odd symplectic
manifolds~\cite{hov:semi},\cite{ass:bv}. We want to emphasize that
\textit{there is no natural volume form preserved by all canonical
transformations}, unlike the even case with the Liouville measure.
However, all ``Darboux coordinate'' volume forms are in the same
orbit of the master groupoid. This orbit is distinguished and it
gives rise to the canonical Laplacian on
\semidensities~\eqref{eqcanlap}, introduced and studied
in~\cite{hov:max,hov:semi}. Every other Laplacian on
half-densities in symplectic case can be expressed via the
canonical operator $\D$ as $\D'=\D-\bro^{-1/2}\D(\bro^{1/2})$.

\textit{For odd symplectic manifolds the modular class vanishes,}
as follows from \eqref{eqdkvadratnafunk}. Every Darboux coordinate
volume form provides a representative of this class by the zero
vector field:  the corresponding Laplacian on functions can be
written as~\eqref{eqcoordeltarho2} and its square evidently
vanishes.  Recall that  Weinstein's modular class on even
symplectic manifolds vanishes  due to Liouville's theorem. The odd
case is more delicate, as there is no analog of the Liouville
volume form. Instead there is a distinguished  class of volume
forms, contained in the same master groupoid orbit.
Lemma~\ref{thmbvlemma} (the ``Batalin--Vilkovisky Lemma'') can be
viewed as a replacement of Liouville's theorem.

We can try to estimate how ``thick''  the orbits of the master
groupoid are and how many such orbits there are. Even in the symplectic
case this analysis is nontrivial.

Consider the orbit through some Darboux coordinate volume form,
i.e., all volume forms $e^{\s(x,\t)}D(x,\t)$ such that
$\D_0\,e^{\s/2}=0$. For infinitesimal $\s$, this reduces
to $\D_0\s=0$. 
Such $\s$ correspond to closed differential forms on a Lagrangian
submanifold~\cite{ass:bv}, \cite{hov:semi} and the dimension of the
orbit is infinite.

On a general odd Poisson manifold $M$, for an arbitrary volume
form $\bro_0$ the modular vector field $\D_{\bro}^2$ is the same
for all points in the orbit of $\bro_0$.  So the orbit of $\bro_0$
is contained in the submanifold $\D_{\bro}^2=\D_{\bro_0}^2$ (fixed
Poisson vector field). It makes sense to  study orbits in such
submanifolds. To estimate their codimension we use
Theorem~\ref{thmzamenakvadrat}: $\DDr-\D_{\bro_0}^2=-
X_{H(\bro,\bro_0)}=0$, which implies that
$H(\bro,\bro_0)=2e^{-{{\s}}/{2}}\D_{\bro_0} e^{{{\s}}/{2}}$ is an
odd Casimir function. Infinitesimally we get that $\D_{\bro_0}\s$
is an odd Casimir. Hence, the codimension of the orbit of $\bro_0$
in the submanifold $\D_{\bro}^2=\const$ equals
\begin{equation}
    \dim\frac{\D_{\bro_0}^{-1}\left\{\text{All odd Casimirs}\right\}}{\D_{\bro_0}^{-1}\{0\}}
    =\dim \left\{\text{All odd Casimirs}\right\}.
\end{equation}

\textit{In the symplectic case}  Casimirs are just constants, and
\textit{the orbits of the master groupoid with the square of the
Laplace operator being fixed  are parametrized by a single
modulus, an odd constant $\nu$}. For example, if $\bro_0$ is a
Darboux coordinate volume form, then
$\nu=\bro^{-1/2}\D\bro^{1/2}$, where $\D$ is the canonical
Laplacian on \semidensities~\eqref{eqcanlap}. This follows from
the analysis performed in other terms in~\cite{hov:semi}.

\subsection{Densities of arbitrary weight}
\label{subsecallweight}

The  commutation formula~\eqref{eqkomdliafun2} for $\D_{\bro}$
acting on functions ($w=0$) is more complicated than the
commutation relation~\eqref{eqkomdliapolu} for \semidensities. The
same is true for the behaviour under a change of  volume form. We
shall now show that this is the generic case and that $w=1/2$ is
an exception.

Consider densities of arbitrary weight $w$. If we fix a basis
volume form $\bro$, then they all have the form
$\bs=s\bro^w$. We define an odd Laplace operator acting on
densities of weight $w$ by the formula
\begin{equation}\label{eqlapw}
    \Dr\bs:=(\Dr s)\bro^w.
\end{equation}
At the r.h.s. stands the Laplacian on functions. The notation here
emphasizes  dependence on $\bro$.

\begin{prop}
On densities of weight $w$ the commutator of the
Laplacian~\eqref{eqlapw} and the multiplication by arbitrary
functions is given by the formula
\begin{equation}\label{eqkomw}
    [\Dr,f]=2\liex{f}+(1-2w)\,\D_{\bro} f.
\end{equation}
\end{prop}

Plugging $w=0$ and $w=1/2$ we recover
formulae~\eqref{eqkomdliafun2} and \eqref{eqkomdliapolu},
respectively. Clearly the case $w=1/2$ is exceptional as the
second term in~\eqref{eqkomw} vanishes identically.

Suppose that we change the basis volume form $\bro$. What is the
transformation law for $\Dr$?

\begin{prop}
On densities of weight $w$ the Laplacian $\Dr$ transforms as
follows:
\begin{equation}
\begin{aligned}\label{eqtransw}
    \D_{\bro'}&=\Dr +(1-2w)\liex{\s}-4w(1-w)\,e^{-\s/2}\Dr
    e^{\s/2}\\
       &=\Dr +(1-2w)\liex{\s}-2w(1-w)\,H(\bro',\bro).
\end{aligned}
\end{equation}
\end{prop}

Plugging $w=0$ and $w=1/2$ we recover formulae~\eqref{eqzamenafun}
and~\eqref{eqzamenarodliapolu1}.

We see again that the  case $w=1/2$ is exceptional, because only
then the transformation~\eqref{eqtransw} involves just a scalar
additive term. In all other cases the transformation includes a
differential operator of the first order. Functions and volume
forms, i.e., $w=0$ and $w=1$, are somewhat special too, because
for them this differential operator  simplifies to a Lie
derivative.

\section{Comparison with Riemannian and even Poisson geometry}

\subsection{Laplacians in  Riemannian geometry}

Let us look at Riemannian geometry from the perspective of odd
Poisson geometry (instead of the converse). Slightly abusing
language, by a Riemannian structure we mean a symmetric tensor
with upper indices $g^{ab}$ not necessarily invertible. The
non-degenerate situation is then an analog of a symplectic
structure. Both Riemannian and odd Poisson structures are
specified by \textit{quadratic Hamiltonians}, even or odd
respectively, which are functions on the same manifold $T^*M$. For
simplicity below  we will write all Riemannian formulae only for
even manifolds, though everything, of course, works in the super
case.

Without assuming that $g^{ab}$ is invertible we can use an
arbitrary volume form $\bro=\rho\,d^nx$ to define a Laplace
operator on functions:
\begin{equation}\label{eqrlap}
    \Dr f:=\divr\grad f=\frac{1}{\r}\,\der{}{x^a}\left(\r g^{ab}\der{f}{x^b}\right).
\end{equation}
The following properties are similar to those of odd
Laplacian~\eqref{eqkomdliafun1}, \eqref{eqlapexp}:
\begin{gather}\label{eqrkomdliafun1}
    \Dr (fg)= (\Dr f)g+f(\Dr g)+2\langle f,g\rangle,\\
    \Dr e^{f}=\Bigl(\Dr f-\langle f,g\rangle\Bigr)\,e^{f}, \label{eqrlapexp}
\end{gather}
where the Poisson bracket is replaced by the ``scalar product of
gradients'' $\langle
f,g\rangle=g^{ab}\,\partial_a{f}\,\partial_b{g}$. Under changes of
volume form the Laplacian~\eqref{eqrlap} transforms as
\begin{equation}\label{eqrzamenafun}
    \Drr = \Dr +\grad {\s},
\end{equation}
where $\brro=e^{\s}\bro$. What fails and has no analogy in the
Riemannian case is the formulae involving the action of the
Laplacian on the bracket.

On densities of arbitrary weight $w$ we can define the Laplacian
by the same formula as~\eqref{eqlapw}:
\begin{equation}\label{eqrlapw}
   \Dr\bs:=(\Dr s)\bro^w,
\end{equation}
if $\bs=s\bro^w$. It has  properties analogous to~\eqref{eqtransw}
and \eqref{eqkomw}.

In particular, the case of \semidensities { }is again
distinguished. Analogs of Theorem~\ref{thmdonsemi} and
Theorem~\ref{thmzamenanapolu} hold. The Laplace operator on
\semidensities { }$\D$ satisfies the condition
\begin{equation}\label{eqrkomsemi}
    [\D,f]=2\lie{\grad {f}}.
\end{equation}
Under the change of volume form
$\bro\to \bro'=e^{\s}\bro$ the Laplacian on \semidensities {
}transforms as
\begin{equation}\label{eqrzamsemi}
    \D'=\D-e^{-{{\s}}/{2}}\Dr e^{{{\s}}/{2}},
\end{equation}
as in~\eqref{eqzamenarodliapolu1}. We again arrive at a groupoid;
the Batalin--Vilkovisky equation is replaced by the Laplace
equation $\D({\bro'}^{1/2})=0$.

Altogether we see that  the analogy between odd Poisson and
``upper'' Riemannian geometry goes unexpectedly far.

\subsection{Geometries controlled by a tensor $T^{ab}$}

Odd Poisson geometry has analogies with Riemannian geometry as
well as with even Poisson geometry.  All three geometries are
controlled by a rank $2$ tensor, say, $T^{ab}$. The difference is
in the type of symmetry and in the parity of $T^{ab}$.

In the even Poisson case, $T^{ab}=P^{ab}$,
$P^{ab}=-(-1)^{\at\bt+\at+\bt}P^{ba}$ and
$\widetilde{P^{ab}}=\at+\bt$. It corresponds to an {even bivector
field} $P\in\fun(\Pi T^*M)$. In the odd Poisson (= Schouten) case,
$T^{ab}=\Sch^{ab}$, $\Sch^{ab}=(-1)^{\at\bt}\Sch^{ba}$, but
$\widetilde{\Sch^{ab}}=\at+\bt+1$. It corresponds to an {odd
quadratic Hamiltonian} $\Sch\in\fun(T^*M)$. In the Riemannian case,
$T^{ab}=g^{ab}$, $g^{ab}=(-1)^{\at\bt}g^{ba}$ and
$\widetilde{g^{ab}}=\at+\bt$. It corresponds to an {even
Hamiltonian} $H=\frac{1}{2}\,g^{ab}p_bp_a\in\fun(T^*M)$.  In all
cases the structure on $M$ is obtained from the  tensor $T$ via
the canonical brackets on $T^*M$ or $\Pi T^*M$.

It is convenient to draw a table:

\par
\bigskip

{\small {\renewcommand{\arraystretch}{1.3} \hspace{-2cm}
\begin{tabular}{|c|c|c|c|} \hline
  \textbf{Even Poisson} &
   \textbf{Odd Poisson (Schouten)}
  & \textbf{Even Riemannian}& \textbf{Odd Riemannian}\\
  $P^{ab}=(-1)^{(\at+1)(\bt+1)}P^{ba}$ & $\Sch^{ab}=(-1)^{\at\bt}\Sch^{ba}$
  & $g^{ab}=(-1)^{\at\bt}g^{ba}$ &
  $\chi^{ab}=(-1)^{(\at+1)(\bt+1)}\chi^{ba}$\\
  \hline
  $P=\frac{1}{2}P^{ab}(x)x^*_bx^*_a$ & $\Sch=\frac{1}{2}\Sch^{ab}(x)p_bp_a$
  & $H=\frac{1}{2}g^{ab}(x)p_bp_a$ &
  $\chi=\frac{1}{2}\chi^{ab}(x)x^*_bx^*_a$\\
  $\Pt=0$ & $\tilde\Sch=1$
  & $\gt=0$ &
  $\tilde\chi=1$ \\
  \parbox{3.6cm}{$\Pi T^*M$,  \small canonical \\
Schouten bracket $\lbr\ ,\ \rbr$} & \parbox{3.4cm}{$T^*M$, \small
canonical \\Poisson bracket $(\ ,\ )$}
  & \parbox{3.4cm}{$T^*M$, \small
canonical \\Poisson bracket $(\ ,\ )$} &
  \parbox{3.6cm}{$\Pi T^*M$,  \small canonical \\
Schouten bracket $\lbr\ ,\ \rbr$}\\[6pt]
\hline \rule{0pt}{22pt}\parbox{3.5cm}{$\{f,g\}=\lbr f,\lbr
P,f\rbr\rbr\\ \vphantom{g^{ab}\der{f}{x^b}}$} &
\rule{0pt}{22pt}\parbox{3.2cm}{$\{f,g\}=(f,(\Sch,g))\\
\vphantom{g^{ab}\der{f}{x^b}}$} &
\rule{0pt}{22pt}\parbox{3.6cm}{$\langle f,g\rangle=(f,(H,g))$\\
\hspace*{0.0cm} $= (-1)^{\ft\at}g^{ab}\der{f}{x^b}\der{g}{x^a}$}
& \rule{0pt}{22pt}\parbox{3.5cm}{$\langle f,g\rangle=\lbr f,\lbr \chi,f\rbr\rbr\\
\vphantom{g^{ab}\der{f}{x^b}}$} \\
$f\mapsto X_f=\{f,\ \}$ & $f\mapsto X_f=(-1)^{\ft+1}\{f,\ \}$ &
$f\mapsto \grad f=\langle f,\ \rangle$ & $f\mapsto
\grad f$ \\
\rule{0pt}{22pt}\parbox{3.5cm}{Jacobi for $\{\ ,\ \}$ \\
$\Leftrightarrow$  $\lbr P,P\rbr=0$} &
\rule{0pt}{22pt}\parbox{3.5cm}{Jacobi for $\{\ ,\ \}$ \\
$\Leftrightarrow$  $(\Sch,\Sch)=0$} & None & None \\[6pt]
\hline
$\Dr f=\divr X_f$ & $\Dr f=\divr X_f$ & $\Dr f=\divr \grad
f$ &
$\Dr f=\divr \grad f$\\
\rule{0pt}{22pt}\parbox{3.5cm}%
{\small $1^{\text{st}}$ order \\
even operator}
& \rule{0pt}{22pt}\parbox{3.5cm}%
{$2^{\text{nd}}$ order \\
odd operator} &
\rule{0pt}{22pt}\parbox{3.5cm}%
{$2^{\text{nd}}$ order \\
even
operator} & \rule{0pt}{22pt}\parbox{3.5cm}{$1^{\text{st}}$ order\\
odd operator}
\\
\hline
\multicolumn{4}{|c|}{$\bro'=e^{\s}\bro$}\\
$\Drr = \Dr +X_{\s}$ & $\Drr = \Dr +X_{\s}$ & $\Drr = \Dr
+\grad{\s}$ & $\Drr = \Dr +\grad{\s}$ \\
 \hline
Modular class: $[\Dr]$ & Modular class: $[\DDr]$ & None & None\\
\hline
\multicolumn{4}{|c|}{Laplacian on half-densities:
$\D(\bs)=\bro^{1/2}\D_{\bro}(\bs\bro^{-1/2})$}\\
\hline nothing good &
\multicolumn{2}{|c|}{\rule{0pt}{22pt}\parbox{5cm}{$\D'-\D=\frac{1}{2}H(\bro',\bro)$\\
$H(\bro',\bro)=e^{-\s/2}\Dr(e^{\s/2})$}} & nothing good \\
\hline None &
\multicolumn{2}{|c|}{\rule{0pt}{22pt}\parbox{7cm}{Master groupoid: \\
$\bro\to\bro'=e^{\s}\bro \ $ such that  $\  \Dr e^{\s/2}=0$}} & None \\
\hline
\end{tabular}
} }

\par\bigskip
The fourth case included in the table but which we shall not develop
here, is that of odd Riemannian structure. Then,
$T^{ab}=\chi^{ab}$ corresponds to an {odd bivector field}. In a
non-degenerate situation, the tensor with lower indices
corresponding to an odd quadratic Hamiltonian $\Sch$ is an odd
$2$-form and that corresponding to an odd bivector field $\chi$ is
an odd symmetric tensor.

Let us comment on some of the similarities and differences.
Algebraic similarity of even/odd Poisson brackets as well as of
even/odd metrics is clear. However, the second order Laplace
operator on functions in the odd Poisson and even Riemannian cases
corresponds to a vector field in the even Poisson and odd
Riemannian cases. Algebraically responsible for arising of the
``master groupoid'' are formulae like $\Dr e^f=e^f\bigl(\Dr
f-\{f,f\}\bigr)$. They do not appear in even Poisson geometry. For
the same reason, in it and in odd Riemannian geometry there is no
interesting theory of Laplace operators acting on various
densities. (Notice that analogies and differences come in pairs
in this picture.)

\subsection{BV formalism and  quantum mechanics}

Without going into details we want to point out  a geometric
analogy between the Batalin--Vilkovisky quantization and ordinary
quantum mechanics.

Let us make a simple but important remark:  the wave function in
the Schr\"o\-dinger equation is a \semidensity. Hence the quantum
Hamiltonian in the Schr\"odinger picture is an operator acting on
\semidensities. Let it have the form $\hat H=-\frac{\hbar^2}{2} \D
+U$, where $\D$ is the Laplace operator on \semidensities. It
satisfies the commutator condition~\eqref{eqrkomsemi}.
$-\frac{\hbar^2}{2} \D$ is a quantization of the   classical
``free'' Hamiltonian $H_0=\frac{1}{2}\,g^{ab}p_bp_a$. The
quasiclassical solution of the Schr\"odinger equation is obtained
by substituting the wave function as
$\boldsymbol{\psi}=e^{iS/\hbar}\bu$ where $\bu=\sum(i\hbar)^n
\bu_n$ is a \semidensity { }and $S$ is a function independent of
$\hbar$. The classical term  is $\langle\grad S,\grad S \rangle$.
The first semiclassical term is $i\hbar\, \lie{\grad S}\bu_0$,
etc. The commutator formula~\eqref{eqrkomsemi} (which is
equivalent to the vanishing of the subprincipal symbol) implies
that the l.h.s. of the transport equations will contain only the
Lie derivative along the gradient of the classical action defined
from the classical Hamilton--Jacobi equation.

Likewise, the ``quantum master equation'' of the
Batalin--Vilkovisky quantization procedure (the
Batalin--Vilkovisky equation) is the equation $\D \bs=0$ for a
\semidensity { }$\bs$ on an odd symplectic manifold, where $\D$ is
the canonical odd Laplacian. Writing formally
$\bs=e^{iS/\hbar}\bu$ as above ($S$ is a function, $\bu$ a
\semidensity) and using the commutator
identity~\eqref{eqkomdliapolu}, one gets a similar expansion in
$i\hbar$ starting from $\{S,S\}=0$ (the ``classical master
equation'').

Hence, we have the following analogy: the quantum master equation
corresponds to the Schr\"odinger equation; its solution (a
\semidensity, often formally written in the purely exponential
form $e^{i\Act/\hbar}$ via the so-called ``quantum effective
action'') corresponds to the wave function; the classical master
equation corresponds to the Hamilton--Jacobi equation or the
eikonal equation; the next quantum corrections, as in the
Schr\"odinger case, are expressed in terms of  Lie derivative
along the Hamiltonian vector field $\liex{S}$, where $S$ is a
solution of the classical master equation.

\smallskip
\textit{Acknowledgement:} The authors want to thank the referee for suggesting numerous 
style improvements of the original manuscript.                

\def\cprime{$'$}

\end{document}